\pdfoutput=1
\documentclass[english,american,reqno]{amsart}
\usepackage[T1]{fontenc}
\usepackage[latin9]{inputenc}
\usepackage{fancyhdr}
\usepackage{refstyle}
\usepackage{units}
\usepackage{textcomp}
\usepackage{amsthm}
\usepackage{amssymb}
\usepackage{setspace}
\usepackage{wasysym}
\usepackage{esint}
\usepackage{prettyref}
\usepackage{xcolor}
\usepackage{babel}
\usepackage{mathrsfs}
\usepackage{url}
\usepackage{amsmath}
\usepackage[numbers]{natbib}
\usepackage[unicode=true,
 bookmarks=true,bookmarksnumbered=false,bookmarksopen=false,
 breaklinks=false,pdfborder={0 0 1},backref=false,colorlinks=true]
 {hyperref}

\makeatletter

\AtBeginDocument{}
\AtBeginDocument{}
\AtBeginDocument{}
\AtBeginDocument{}
\AtBeginDocument{}
\AtBeginDocument{\providecommand\ref[1]{\ref{eq:#1}}}
\AtBeginDocument{}
\RS@ifundefined{subref}
  {\def\RSsubtxt{section~}\newref{sub}{name = \RSsubtxt}}
  {}
\RS@ifundefined{thmref}
  {\def\RSthmtxt{theorem~}\newref{thm}{name = \RSthmtxt}}
  {}
\RS@ifundefined{lemref}
  {\def\RSlemtxt{lemma~}\newref{lem}{name = \RSlemtxt}}
  {}

\numberwithin{equation}{section}
\numberwithin{figure}{section}
\theoremstyle{plain}
\newtheorem{theorem}{\protect\theoremname}[section]
  \theoremstyle{plain}
  \newtheorem{proposition}[theorem]{\protect\propositionname}
  \theoremstyle{plain}
  \newtheorem{lemma}[theorem]{\protect\lemmaname}
  \theoremstyle{remark}
  \newtheorem{remark}[theorem]{\protect\remarkname}
  \theoremstyle{plain}
  \newtheorem{corollary}[theorem]{\protect\corollaryname}

\makeatother

\usepackage{babel}
  \addto\captionsamerican{\renewcommand{\corollaryname}{Corollary}}
  \addto\captionsamerican{\renewcommand{\lemmaname}{Lemma}}
  \addto\captionsamerican{\renewcommand{\propositionname}{Proposition}}
  \addto\captionsamerican{\renewcommand{\remarkname}{Remark}}
  \addto\captionsamerican{\renewcommand{\theoremname}{Theorem}}
  \addto\captionsenglish{\renewcommand{\corollaryname}{Corollary}}
  \addto\captionsenglish{\renewcommand{\lemmaname}{Lemma}}
  \addto\captionsenglish{\renewcommand{\propositionname}{Proposition}}
  \addto\captionsenglish{\renewcommand{\remarkname}{Remark}}
  \addto\captionsenglish{\renewcommand{\theoremname}{Theorem}}
  \providecommand{\corollaryname}{Corollary}
  \providecommand{\lemmaname}{Lemma}
  \providecommand{\propositionname}{Proposition}
  \providecommand{\remarkname}{Remark}
\providecommand{\theoremname}{Theorem}

\title{Spaces of harmonic maps of the projective plane to the four-dimensional
sphere}

\author{Ravil Gabdurakhmanov}
\address{Faculty of Mathematics, National Research University
Higher School of Economics, 6 Usacheva Str., 119048, Moscow, Russia }
\address{Independent University of Moscow, Bolshoy Vlasyevskiy
Pereulok 11, 119002, Moscow, Russia}
\email{ravil.gabdurakhmanov@gmail.com}
\curraddr{School of Mathematics, University of Leeds, Leeds, LS2 9JT, United Kingdom}
\thanks{Supported in part by the Simons Foundation.}
\begin{document}
\begin{abstract}
The spaces of harmonic maps of the projective plane to the four-dimensional
sphere are investigated in this paper by means of twistor lifts. It
is shown that such spaces are empty in case of even harmonic degree.
In case of harmonic degree less than 6 it was shown that such spaces
are path-connected and an explicit parameterization of the canonical
representatives was found. In addition, the last section summarizes
known results for harmonic maps of the two-dimensional sphere to the
four-dimensional sphere of harmonic degree less than 6.
\end{abstract}

\maketitle

\section{Introduction}

Let $\phi:M\rightarrow N$ be a map between Riemannian manifolds.
We define its \textit{energy} by the formula 
\[
E(\phi)=\frac{1}{2}\intop_{M}|d\phi(x)|^{2}dx,
\]

\noindent where $d\phi(x)$ is the differential of $\phi$ at the
point $x\in M$; and $dx$ is the volume element of $M$. Euler-Lagrange
operator $\tau(\phi)=div(d\phi)$ associated with the functional $E$
is called a \textit{tension field of} $\phi$. The map $\phi:M\rightarrow N$
is said to be \textit{harmonic} if its tension field vanishes identically
i.e. $\phi$ is a critical point of $E$.

Harmonic maps appear in mathematics in different ways, in the theory
of minimal surfaces, in spectral geometry, where theory of harmonic
maps from $S^{2}$ to $S^{n}$ led to major advances in isoperimetric
inequalities for eigenvalues of the Laplace-Beltrami operator, see
e.g. \citep{key-14},\citep{key-15}, etc. For example, there are
the following particular cases
\begin{itemize}
\item If $dim\,M=1$, then the harmonic maps are the geodesics of $N$.
\item If $N=\mathbb{R}$, they are harmonic functions on $M$.
\item If $dim\,M=2$, they include parametric representations of the minimal
surfaces of $N$; the energy is the Dirichlet-Douglas integral.
\end{itemize}
As a result, there are many papers on the theory of the harmonic maps.
This paper is based on the results of the famous works of E. Calabi
\citep{key-1,key-2}, J. Barbosa \citep{key-3}, R. L. Bryant \citep{key-4},
and essentially uses ideas developed in works of J. Bolton and L.
M. Woodward \citep{key-6,key-7,key-8}. We use the following proposition
and fundamental theorems.
\begin{proposition}
\textup{\label{prop:-An-isometric}\citep{key-9} }An isometric immersion
$\phi:(M,g)\rightarrow(N,h)$ is minimal if and only if it is harmonic.
\end{proposition}
\begin{theorem}
\textup{\label{thm:(E.-Calabi-).}(}E. Calabi\textup{ \citep{key-1}).}
Let $X:S^{2}\rightarrow\mathbb{R}^{n}$ be an immersion of a sphere
$S^{2}$ into $\mathbb{R}^{n}$, whose image is a locally minimal
surface in a sphere $rS^{n-1}$ of radius $r$, and is not contained
in any hyperplane of $\mathbb{R}^{n}$. Then the following conclusions
hold.

(1) The area $A=A(S^{2})$ of the image sphere $S^{2}$ is an integer
multiple of $2\pi r^{2}$.

(2) The dimension $n$ is odd.
\end{theorem}
In this theorem one uses an assumption that the image of $X:S^{2}\rightarrow\mathbb{R}^{n}$
is not contained in any hyperplane of $\mathbb{R}^{n}$. Such maps
are called \textit{linearly full. }It should be noted that if we induce
metric on $S^{2}$ by \foreignlanguage{english}{$X$} in the conditions
of \prettyref{thm:(E.-Calabi-).}, then $X$ is an isometric immersion,
and hence harmonic due to Proposition \prettyref{prop:-An-isometric}.
J. Barbosa clarifies paragraph (1) of \prettyref{thm:(E.-Calabi-).}
in the case of $r=1$.
\begin{theorem}
\textup{\label{(J.-Barbosa-).}(}J. Barbosa\textup{ \citep{key-3}).}
The area of a generalized minimal immersion $X:S^{2}\rightarrow S^{2m}$
is $4\pi d$, for some $d\in\mathbb{N}$, $d\geq\frac{m(m+1)}{2}$.
\end{theorem}
The word \textit{generalized} in this theorem means that the immersion
$X$ may be branched at isolated points. The number $d$, appearing
in this theorem, is called a \textit{harmonic degree} (or \textit{degree})
of a harmonic map $X$. The case $m=2$ has been investigated by R.
L. Bryant in \citep{key-4} using the construction of twistor fibration
$\pi:\,\mathbb{C}P^{3}\rightarrow S^{4}$ introduced in Section 2.
In the following theorem \textsl{horizontal }means orthogonal to the
fibers of $\pi$ (see \prettyref{sec:2}).
\begin{theorem}
\textup{\label{thm:(R.-L.-Bryant}(}R. L. Bryant\textup{ \citep{key-4})}. 

(1) Let $\psi:S^{2}\rightarrow\mathbb{C}P^{3}$ be a linearly full
horizontal holomorphic curve then $\phi=\pi\psi:S^{2}\rightarrow S^{4}$
is a linearly full harmonic map. Conversely every linearly full harmonic
map $\phi:S^{2}\rightarrow S^{4}$ is equal to $\pm\pi\psi$ for some
uniquely determined linearly full horizontal holomorphic curve $\psi$
called the \textup{twistor lift} of $\phi$. 

(2) Twistor lift of a linearly full harmonic map $\phi:S^{2}\rightarrow S^{4}$
is an algebraic curve.

(3) Area of $S^{2}$ w.r.t. the metric induced by $\phi$ is equal
to $4\pi d$, and $d$ is a degree of its twistor lift.
\end{theorem}
This paper investigates harmonic maps $\phi:\mathbb{R}P^{2}\rightarrow S{}^{4}$
by using standard Riemannian double cover $S^{2}\rightarrow\mathbb{R}P^{2}$
and twistor lifts to $\mathbb{C}P^{3}$ of \textit{invariant} harmonic
maps $S^{2}\rightarrow S^{4}$ (i.e. harmonic maps that can be factorized
through the double cover $S^{2}\rightarrow\mathbb{R}P^{2}$). In this
way we define harmonic degree of the map $\phi:\mathbb{R}P^{2}\rightarrow S{}^{4}$
as a degree of its lift $\tilde{\phi}:S^{2}\rightarrow S{}^{4}$ under
double cover.

The structure of the paper is following. In \prettyref{sec:2} we
give an overview of the twistor fibration and a groups acting on it.
\prettyref{sec:3} is devoted to a higher singularities that occur
for the twistor lifts of harmonic maps. In \prettyref{sec:4} topologies
on the sets of linearly full horizontal holomorphic curves and harmonic
maps are defined, and their homeomorphism is described. \prettyref{sec:5}
addresses general results on harmonic maps of arbitrary degree, for
example, it shows that there are no harmonic maps from the projective
plane to the four-dimensional sphere of even degree. In \prettyref{sec:6}
the new facts obtained in previous section and similar ideas as used
by J. Bolton and L. M. Woodward in \citep{key-6} allow us to describe
the spaces of harmonic maps of degree less than 6, and to show path-connectedness
of these spaces and the ``bubbling'' phenomenon. \prettyref{sec:7}
summarizes some known results for harmonic maps from the two-dimensional
sphere to the four-dimensional sphere.

\section{\label{sec:2}The twistor fibration}

We introduce the construction of the \textit{twistor fibration} $\pi:\,\mathbb{C}P^{3}\rightarrow S^{4}$.
The \textit{Hopf map} $\rho:\,\mathbb{C}P^{3}\rightarrow\mathbb{H}P^{1}$
is given by
\[
\rho([z_{1},z_{2},z_{3},z_{4}])=[z_{1}+z_{2}j,z_{3}+z_{4}j].
\]

\noindent The canonical identification of $\mathbb{H}P^{1}$ and $S^{4}\subset\mathbb{H}\oplus\mathbb{R}=\mathbb{R}^{5}$
is given by the stereographic projection from the south pole of $S^{4}$
onto the equatorial hyperplane $\mathbb{H}$ in $\mathbb{R}^{5}$
which is included in $\mathbb{H}P^{1}$ by $q\mapsto[1,q]$. This
identification is given by
\[
[q_{1},q_{2}]\in\mathbb{H}P^{1}\leftrightarrow\frac{(2\bar{q}_{1}q_{2},|q_{1}|^{2}-|q_{2}|^{2})}{|q_{1}|^{2}+|q_{2}|^{2}}\in S^{4}.
\]

\noindent Now if $\mathbb{H}^{2}$ is a left quaternionic vector space,
then $\pi$ is obtained by composing Hopf map $\rho$ with the canonical
identification of $\mathbb{H}P^{1}$ and $S^{4}$.

\noindent Hence, we obtain an explicit formula for the twistor fibration
\[
\pi([z_{1},z_{2},z_{3},z_{4}])=\frac{(2(\bar{z}_{1}z_{3}+z_{2}\bar{z}_{4}),2(\bar{z}_{1}z_{4}-z_{2}\bar{z}_{3}),|z_{1}|^{2}+|z_{2}|^{2}-|z_{3}|^{2}-|z_{4}|^{2})}{|z_{1}|^{2}+|z_{2}|^{2}+|z_{3}|^{2}+|z_{4}|^{2}}.
\]

If $\mathbb{C}P^{3}$ is endowed with the Fubini-Study metric of constant
holomorphic sectional curvature 1 then $\pi$ is a Riemannian submersion
and the \textit{horizontal distribution} on $\mathbb{C}P^{3}$ consists
of those tangent vectors to $\mathbb{C}P^{3}$ which are orthogonal
to the fibers of $\pi$.

Let 
\begin{equation}
J=\begin{pmatrix}0 & -1 & 0 & 0\\
1 & 0 & 0 & 0\\
0 & 0 & 0 & -1\\
0 & 0 & 1 & 0
\end{pmatrix}.
\end{equation}

\noindent Then the projectivization $PSp(2;\mathbb{C})$ of the complexified
symplectic group 
\[
Sp(2;\mathbb{C})=\{A\in SL(4;\mathbb{C})|A^{t}JA=J\}
\]

\noindent acts on $\mathbb{C}P^{3}$ in the standard way as the group
of holomorphic diffeomorphisms, which preserve the horizontal distribution,
with $PSp(2)=PSp(2;\mathbb{C})\cap PU(4)$ the subgroup of holomorphic
isometries which preserve the horizontal distribution \citep{key-7}. 

It should be noted that all of these groups act transitively on $\mathbb{C}P^{3}$
and $Sp(2)$ and hence $PSp(2)$ are path-connected \citep{key-10}.
It should be noted also that the elements of $Sp(2)$ are matrices
$(\mathbf{u},J\mathbf{\bar{u}},\mathbf{v},J\mathbf{\bar{v}})$ such
that $\mathbf{u},\mathbf{v}$ are unit vectors in $\mathbb{C}^{4}$,
and $\mathbf{v}$ unitarily orthogonal to $\mathbf{u}$ and $J\mathbf{\bar{u}}$.

Let $S$ be a Riemann surface and $\psi:S\rightarrow\mathbb{C}P^{3}$
be a holomorphic curve. Then $\psi$ is \textit{horizontal} if at
each point it is tangent to the horizontal distribution on $\mathbb{C}P^{3}$,
or equivalently if it intersects each fiber of $\pi$ orthogonally.
It was proved in \citep{key-4} that the condition for a holomorphic
curve $\psi=[\mathbf{f}]=[f_{1},f_{2},f_{3},f_{4}]$ to be horizontal
is 
\begin{equation}
(\mathbf{f}^{\prime},J\mathbf{f})=0,
\end{equation}

\noindent where $(\,,\,)$ denotes the complex bilinear extension
to $\mathbb{C}^{4}$ of the standard real scalar product on $\mathbb{R}^{4}$.
The horizontality condition may be written as
\begin{equation}
f_{1}^{\prime}f_{2}-f_{1}f_{2}^{\prime}+f_{3}^{\prime}f_{4}-f_{3}f_{4}^{\prime}=0.\label{eq:2.3}
\end{equation}

\noindent Here and throughout the rest of the paper we write holomorphic
curves $\psi:S\rightarrow\mathbb{C}P^{n}$ locally in terms of a local
complex coordinate $z$ in the form $\psi=[\mathbf{f}]=[f_{1},f_{2},...,f_{n+1}]$,
where $f_{1},...,f_{n+1}$ are holomorphic functions of $z$ without
common zeros. We will call such a representation a \textit{reduced
form} of $\psi$.

\section{\label{sec:3}Higher singularities for horizontal holomorphic curves}

Let us recall the definition of singularity type of holomorphic curves
in $\mathbb{C}P^{n}$ \citep{key-6}. Let $S$ be a Riemann surface
and suppose that $\psi:S\rightarrow\mathbb{C}P^{n}$ is a linearly
full holomorphic curve. We write $\psi$ locally in reduced form $\psi=[\mathbf{f}]=[f_{1},f_{2},...,f_{n+1}]$,
and let $\mathbf{f}^{(j)}$ denote the $j$-th derivative of $\mathbf{f}$
with respect to $z$. For each $k=0,\ldots,n-1$ we define $k$\textit{-th
osculating curve} as $[\mathbf{f}\land\ldots\land\mathbf{f}^{(k)}]$.
A \textit{higher singularity} of $\psi$ is a point $p\in S$ where
the derivative of $k$-th osculating curve is equal to zero for some
$k=0,\ldots,n-1$. The set $Z(\psi)$ of higher singularities of $\psi$
is therefore given by
\[
Z(\psi)=\{p\in S\,|\,\mathbf{f}\land\ldots\land\mathbf{f}^{(n)}(p)=0\}.
\]

\noindent If $z$ is a local complex coordinate with $z(p)=0$, then
$\mathbf{f}$ can be written in the \textit{normal form}
\[
\mathbf{f}(z)=h_{0}(z)\mathbf{a}_{0}+z^{r_{0}(p)+1}h_{1}(z)\mathbf{a}_{1}+\ldots+z^{r_{0}(p)+\ldots r_{n-1}(p)+n}h_{n}(z)\mathbf{a}_{n}
\]

\noindent for some suitable choice of basis $\mathbf{a}_{0},\ldots,\mathbf{a}_{n}$
of $\mathbb{C}^{n+1}$, non-negative integers\linebreak{}
 $r_{0}(p),\ldots,r_{n-1}(p)$, and holomorphic functions $h_{0}(z),\ldots,h_{n}(z)$,
each non-zero at $z=0$. If $r_{k}(p)>0$ then the derivative of the
$k$-th osculating curve has a zero of order $r_{k}(p)$ at $p$ and
the \textit{singularity type} of $\psi$ is defined to be the set
\[
\{(p;r_{0}(p),\ldots,r_{n-1}(p))\,|\,p\in Z(\psi)\}.
\]

\begin{lemma}
\label{lem:-Let-}\textup{\citep{key-6}} Let $g\in PGL(n+1;\mathbb{C})$
and $\omega$ be a conformal diffeomorphism of $S$. If $\psi$ has
singularity type $\{(p;r_{0}(p),\ldots,r_{n-1}(p))\,|\,p\in Z(\psi)\}$
then $g\psi\omega^{-1}$ has singularity type $\{(\omega(p);r_{0}(p),\ldots,r_{n-1}(p))\,|\,p\in Z(\psi)\}$.
\end{lemma}
We now give a criterion for determining the higher singularities of
a linearly full horizontal holomorphic curve $\psi:S\rightarrow\mathbb{C}P^{3}$.
\begin{lemma}
\textup{\citep{key-6}} Let $\psi:S\rightarrow\mathbb{C}P^{3}$ be
a linearly full horizontal holomorphic curve written locally in terms
of a complex coordinate $z$ in the reduced form as $\psi(z)=[\mathbf{f}(z)]=[f_{1}(z),f_{2}(z),f_{3}(z),f_{4}(z)]$.
Then $z=a$ is a higher singularity of $\psi$ if and only if
\begin{equation}
(f_{1}^{\prime\prime}f_{2}^{\prime}-f_{1}^{\prime}f_{2}^{\prime\prime}+f_{3}^{\prime\prime}f_{4}^{\prime}-f_{3}^{\prime}f_{4}^{\prime\prime})(a)=0,\label{eq:3.1}
\end{equation}

\noindent or equivalently, $(\mathbf{f}^{\prime\prime},J\mathbf{f}^{\prime})(a)=0$.
\end{lemma}
For a linearly full holomorphic curve $\psi:S^{2}\rightarrow\mathbb{C}P^{n}$
and $k=0,\ldots,n-1$ we define 
\[
r_{k}=\sum_{p\in Z(\psi)}r_{k}(p).
\]

\begin{proposition}
Let $\psi:S^{2}\rightarrow\mathbb{C}P^{3}$ be a linearly full horizontal
holomorphic curve and $p\in S^{2}$. Then $r_{0}(p)=r_{2}(p)$ and
\begin{equation}
2r_{0}+r_{1}=2d-6,\label{eq:3.2}
\end{equation}

\noindent where $d$ is the degree of $\psi$.
\end{proposition}
Proofs of the lemmas and proposition of this section can be found
in \citep{key-6}.

\section{\label{sec:4}A canonical form}

Due to \prettyref{thm:(R.-L.-Bryant}, the fact that linearly full
harmonic map $\phi:S^{2}\rightarrow S^{4}$ has a twistor lift which
is, moreover, an algebraic curve, gives us a good opportunity for
investigating harmonic maps. Namely, let $HHol_{d}^{LF}(\mathbb{C}P^{3})$
denote the set of linearly full horizontal holomorphic maps of $S^{2}$
to $\mathbb{C}P^{3}$ of degree $d$, and let $Harm_{d}^{LF}(S^{4})$
denote the set of linearly full harmonic maps of $S^{2}$ to $S^{4}$
with induced area $4\pi d$. Then $Harm_{d}^{LF}(S^{4})$ is divided
into two components $Harm_{d}^{+}(S^{4})$ and $Harm_{d}^{-}(S^{4})$,
which can be interchanged by the antipodal involution of $S^{4}$,
and there is a \foreignlanguage{british}{bijective} correspondence

\begin{equation}
\pi_{*}^{\pm}:HHol_{d}^{LF}(\mathbb{C}P^{3})\rightarrow Harm_{d}^{\pm}(S^{4})
\end{equation}

\noindent defined by $\pi_{*}^{\pm}(\psi)=\pm\pi\circ\psi$. 

Clearly, $p\in S^{2}$ is a branch point of $\psi$ if and only if
$r_{0}(p)>0$ while as is shown in Section 7 of \citep{key-7}, if
$r_{0}(p)=0$ then $p$ is an umbilic of $\psi$ if and only if $r_{1}(p)>0$.
Further, it follows from the paper \citep{key-11} that the branch
points and umbilics of $\phi$ coincide with those of $\psi$. Thus
the higher singularities of $\psi$ occur at the branch points or
umbilics of $\phi$.

It is clear that the group $PSp(2)$ acts freely on $HHol_{d}^{LF}(\mathbb{C}P^{3})$
by postcomposition via its standard action on $\mathbb{C}P^{3}$ and
it follows from \prettyref{lem:-Let-} that the singularity type is
preserved. The identifications determined by $\pi_{*}^{\pm}$ determine
an action of $PSp(2)$ on $Harm_{d}^{\pm}(S^{4})$ which preserves
induced area, branch points, umbilics and antipodal invariance. Similarly
the rotation group $SO(3)\cong PSU(2)$ also acts on $HHol_{d}^{LF}(\mathbb{C}P^{3})$
and $Harm_{d}^{\pm}(S^{4})$ by precomposition via its standard action
on $S^{2}=\mathbb{C}P^{1}$ and also preserves antipodal invariance.
Note that the actions of $PSp(2)$ and $PSU(2)$ commute and that
the maps $\pi_{*}^{\pm}$ are $PSp(2)$-equivariant and $PSU(2)$-equivariant.

Now we consider the vector space $\mathbb{C}[z]_{d}$ of polynomials
of degree at most $d$, and let $V$ be the subset of $(\mathbb{C}[z]_{d})^{4}$
consisting of those quadruples of coprime polynomials with maximum
degree $d$ for which the map $z\rightarrow[f_{1}(z),f_{2}(z),f_{3}(z),f_{4}(z)]$
is linearly full in $\mathbb{C}P^{3}$. Then $V$ is a projective
subset of $(\mathbb{C}[z]_{d})^{4}$, and we identify its projectivisation
$P(V)$ with the space of linearly full holomorphic maps of degree
$d$ from $S^{2}$ to $\mathbb{C}P^{3}$ in the usual way via
\[
[f_{1},f_{2},f_{3},f_{4}]\leftrightarrow(z\rightarrow[f_{1}(z),f_{2}(z),f_{3}(z),f_{4}(z)]).
\]

\noindent Here, and subsequently, we use the complex coordinate on
$S^{2}$ defined by the stereographic projection from the south pole
of $S^{2}$ onto the equatorial plane which means, in the usual sense,
we may identify $S^{2}$ with $\mathbb{C}\cup\{\infty\}$. According
to this identification antipodal involution interchanges $z$ and
$-\frac{1}{\bar{z}}$.

We endow $(\mathbb{C}[z]_{d})^{4}$ with its natural topology as a
vector space, and $P((\mathbb{C}[z]_{d})^{4})$ the quotient topology.
Then $V$ is an open subset of $(\mathbb{C}[z]_{d})^{4}$, and $P(V)$
is an open subset of $P((\mathbb{C}[z]_{d})^{4})$. Subsets of any
of these spaces are then endowed with the induced topology. We denote
by $V^{H}\subset V$ the subset of horizontal maps. Next we give $Harm_{d}^{LF}(S^{4})$
the compact-open topology and rewrite a lemma from \citep{key-8}.
\begin{lemma}
\label{lem:4.1}$Harm_{d}^{\pm}(S^{4})$ are closed subsets of $Harm_{d}^{LF}(S^{4})$,
and the maps $\pi_{*}^{\pm}:P(V^{H})\rightarrow Harm_{d}^{\pm}(S^{4})$
are homeomorphisms.
\end{lemma}
\noindent For convenience, we introduce \textit{coefficient matrix
}(or\textit{ matrix}) of a holomorphic curve $\mathbf{f}(z)=(f_{1}(z),f_{2}(z),f_{3}(z),f_{4}(z))\in(\mathbb{C}[z]_{d})^{4}$.
It is given by
\[
F=\begin{pmatrix}a_{1}^{0} & a_{1}^{1} & \ldots & a_{1}^{d-1} & a_{1}^{d}\\
a_{2}^{0} & a_{2}^{1} & \ldots & a_{2}^{d-1} & a_{2}^{d}\\
a_{3}^{0} & a_{3}^{1} & \ldots & a_{3}^{d-1} & a_{3}^{d}\\
a_{4}^{0} & a_{4}^{1} & \ldots & a_{4}^{d-1} & a_{4}^{d}
\end{pmatrix},
\]

\noindent where $a_{i}^{j}$ is a $j$-th coefficient of the polynomial
$f_{i}(z)=a_{i}^{0}+a_{i}^{1}z+\ldots+a_{i}^{j}z^{j}+\ldots+a_{i}^{d}z^{d}$,
and $1\leq i\leq4$, $0\leq j\leq d$. It is clear from definition
that $\mathbf{(f}(z))^{t}=F\cdot(1,z,z^{2},\ldots,z^{d})^{t}$ and
also that $\mathbf{f}(z)$ is linearly full if the coefficient matrix
$F$ has a full rank.

Note that an element $A\in Sp(2)$, acting on $\mathbf{f}(z)\in(\mathbb{C}[z]_{d})^{4}$
by postcomposition (i.e. $\mathbf{f}(z)\rightarrow A\circ\mathbf{f}(z)$),
acts on coefficient matrix $F$ by the left multiplication (i.e. $F\rightarrow AF$).
Also, as noted earlier, the group $PSU(2)\cong SO(3)$ acts on $S^{2}$
as a group of rotations. This action defines an action of $PSU(2)$
on $\mathbb{R}P^{2}$ since it preserves the antipodal points. It
should be noted that the group $PSU(2)$ is path-connected \citep{key-10}. 

In the following sections we investigate linearly full harmonic maps
of $S^{2}$ to $S^{4}$ of degrees $d=3$, $4$, $5$ and $6$, which
are the lifts of harmonic maps of $\mathbb{R}P^{2}$ to $S^{4}$ through
the double cover, i.e. maps invariant w.r.t. the antipodal map (we
will call them \textsl{invariant maps }for brevity). Let $AHarm_{d}^{LF}(S^{4})\subset Harm_{d}^{LF}(S^{4})$
denote the space of invariant linearly full harmonic maps of $S^{2}$
to $S^{4}$of degree $d$, with $AHarm_{d}^{\pm}(S^{4})\subset Harm_{d}^{\pm}(S^{4})$.
And let us denote the twistor lifts of this spaces as $AHHol_{d}^{LF}(\mathbb{C}P^{3})\subset HHol_{d}^{LF}(\mathbb{C}P^{3})$.
Our main idea is the deformation of elements of $AHarm_{d}^{\pm}(S^{4})\subset Harm_{d}^{\pm}(S^{4})$
to some canonical form by the action by some appropriate elements
of $PSU(2)$ and $PSp(2)$. Note that these groups act continuously
on $AHHol_{d}^{LF}(\mathbb{C}P^{3})$ and, due to their path-connectedness,
we can recover some information about path-connectedness of $AHHol_{d}^{LF}(\mathbb{C}P^{3})$.

\section{\label{sec:5}Harmonic maps of arbitrary degree}

In this section we will exclude the condition of antipodal invariance
by finding the general form of invariant maps. In this way we will
need some addition to \prettyref{thm:(R.-L.-Bryant}.
\begin{theorem}
\textup{\label{thm:(R.-L.-Bryant-1}(}R. L. Bryant\textup{ \citep{key-4})}. 

Let $\hat{\psi}:S^{2}\rightarrow\mathbb{C}P^{3}$ be a linearly full
horizontal antiholomorphic curve then $\phi=\pi\hat{\psi}:S^{2}\rightarrow S^{4}$
is a linearly full harmonic map. Conversely every linearly full harmonic
map $\phi:S^{2}\rightarrow S^{4}$ is equal to $\pm\pi\hat{\psi}$
for some uniquely determined linearly full horizontal antiholomorphic
curve $\hat{\psi}$ called the \textup{antiholomorphic} \textup{twistor
lift} of $\phi$.

Moreover if $\psi=[\mathbf{f}]=[f_{1},f_{2},f_{3},f_{4}]$ is a twistor
lift of $\phi$ then $\hat{\psi}=\bar{\psi}J=[\bar{\mathbf{f}}J]=[\bar{f_{2}},-\bar{f_{1}},\bar{f_{4}},-\bar{f_{3}}]$
is an antiholomorphic twistor lift of $\phi$ \citep{key-7}.
\end{theorem}
For $\mathbf{f}(z)=(f_{1}(z),f_{2}(z),f_{3}(z),f_{4}(z))$ of degree
$n$ we introduce 
\[
\begin{array}{c}
\tilde{\mathbf{f}}(z)=(\tilde{f_{1}}(z),\tilde{f_{2}}(z),\tilde{f_{3}}(z),\tilde{f_{4}}(z))=\\
=(\bar{z}^{n}f_{1}(-1/\bar{z}),\bar{z}^{n}f_{2}(-1/\bar{z}),\bar{z}^{n}f_{3}(-1/\bar{z}),\bar{z}^{n}f_{4}(-1/\bar{z})).
\end{array}
\]
 Then $[\tilde{f}(z)]$ is an antiholomorphic curve in $\mathbb{C}P^{3}$.
\begin{theorem}
\label{thm:5.2}There are no harmonic maps of real projective plane
to the four-dimensional sphere of even degree, i.e. spaces $AHarm_{n}^{LF}(S^{4})$
are empty for even $n$.
\end{theorem}
\begin{proof}
Let us suppose that such maps exist. Then for each map $\phi\in AHarm_{n}^{LF}(S^{4})$
we will compose its twistor lift $\psi=[\mathbf{f}(z)]=[f_{1}(z),f_{2}(z),f_{3}(z),f_{4}(z)]$
with antipodal involution of $S^{2}$ and get an antiholomorphic twistor
lift $\tilde{\psi}=[\mathbf{\tilde{f}}(z)]$. Using \prettyref{thm:(R.-L.-Bryant-1}
we get $\psi=-\bar{\tilde{\psi}}J$. Then 
\begin{equation}
\begin{array}{c}
g(z)\cdot\mathbf{f}(z)=g(z)\cdot(f_{1}(z),f_{2}(z),f_{3}(z),f_{4}(z))=\\
=-\overline{\tilde{\mathbf{f}}(z)}J=(-\overline{\tilde{f_{2}}(z)},\overline{\tilde{f_{1}}(z)},-\overline{\tilde{f_{4}}(z)},\overline{\tilde{f_{3}}(z)}),
\end{array}\label{eq:5.1}
\end{equation}
 for some function $g(z)$ on $S^{2}$. Now we will prove that $g(z)$
is a holomorphic function. From \ref{eq:5.1} we have

\begin{equation}
g(z)=\frac{-\overline{\tilde{f_{2}}(z)}}{f_{1}(z)}=\frac{\overline{\tilde{f_{1}}(z)}}{f_{2}(z)}=\frac{-\overline{\tilde{f_{4}}(z)}}{f_{3}(z)}=\frac{\overline{\tilde{f_{3}}(z)}}{f_{4}(z)}.\label{eq:5.2}
\end{equation}
It is clear that $\overline{\tilde{f_{i}}(z)}$ $(i=1,2,3,4)$ are
polynomials of $z$ of degree less or equal to $n$, so they are holomorphic
functions on $\mathbb{C}$. At the same time $f_{i}(z)$ $(i=1,2,3,4)$
don't have common roots and therefore $g(z)$ is holomorphic on $\mathbb{C}$.
By the definition there is polynomial of degree $n$ among $f_{i}(z)$
$(i=1,2,3,4)$ and from \ref{eq:5.2} we conclude that $g(z)$ is
holomorphic at $\infty$. Thus $g(z)$ is the holomorphic function
on $S^{2}$ and therefore $g(z)=\frac{1}{\beta}=const\neq0$.

Since $f_{i}(z)=a_{i}^{0}+a_{i}^{1}z+\ldots+a_{i}^{j}z^{j}+\ldots+a_{i}^{n}z^{n}$,
$(i=1,2,3,4)$ we have $\overline{\tilde{f_{i}}(z)}=(-1)^{n}\overline{a_{i}^{n}}+(-1)^{n-1}\overline{a_{i}^{n-1}}z+\ldots+(-1)^{n-j}\overline{a_{i}^{n-j}}z^{j}+\ldots+\overline{a_{i}^{0}}z^{n}$,
$(i=1,2,3,4)$. Considering \ref{eq:5.2} we can write

\begin{equation}
\begin{cases}
a_{1}^{n}=-\beta\overline{a_{2}^{0}}, & a_{2}^{0}=(-1)^{n}\beta\overline{a_{1}^{n}},\\
a_{2}^{n}=\beta\overline{a_{1}^{0}}, & a_{1}^{0}=(-1)^{n+1}\beta\overline{a_{2}^{n}},\\
a_{3}^{n}=-\beta\overline{a_{4}^{0}}, & a_{4}^{0}=(-1)^{n}\beta\overline{a_{3}^{n}},\\
a_{4}^{n}=\beta\overline{a_{3}^{0}}, & a_{3}^{0}=(-1)^{n+1}\beta\overline{a_{4}^{n}}.
\end{cases}
\end{equation}

From these relations we can conclude for leading coefficients $a_{i}^{n}=(-1)^{n+1}|\beta|^{2}a_{i}^{n}$
and since $n$ is even one has $a_{i}^{n}=-|\beta|^{2}a_{i}^{n}$.
But it can only be when $a_{i}^{n}=0$ $(i=1,2,3,4)$ and thus the
degree of the curve is less than n. So we have a contradiction which
proves this theorem.
\end{proof}

\begin{remark}
\label{rem:In-the-case}In case of odd $n$ we have $a_{i}^{n}=|\beta|^{2}a_{i}^{n}$
for leading coefficients, therefore $|\beta|=1$ and $\beta=\mathrm{exp}(i\tilde{\beta})$
for some $\tilde{\beta}\in\mathbb{R}$. In this way any curve of odd
degree is of the form $\psi=[\mathbf{f}(z)]=[f_{1}(z),\beta\overline{\tilde{f_{1}}(z)},f_{3}(z),\beta\overline{\tilde{f_{3}}(z)}]$.
Note that if we act on the curve by rotation $\omega(z)=\mathrm{exp(-i\frac{\tilde{\beta}}{n}})z$,
we get the curve with $\beta=1$. (Note also that this rotation leaves
points $0$ and $\infty$ unchanged). Thus any curve of odd degree
can be rotated to a curve of form $\psi=[\mathbf{f}(z)]=[f_{1}(z),\overline{\tilde{f_{1}}(z)},f_{3}(z),\overline{\tilde{f_{3}}(z)}]$.

We need following useful lemma in next sections. 
\end{remark}
\begin{lemma}
\label{lem:5.4}Any curve $\psi\in HHol_{d}^{LF}(\mathbb{C}P^{3})$
by an appropriate element $g\in PSp(2)\subset PSp(2,\mathbb{C})$
can be reduced to the curve $g\psi$ with coefficient matrix of form
\[
\begin{pmatrix}\hexstar & \hexstar & \hexstar & \hexstar & ... & \hexstar & \hexstar & \hexstar & \hexstar\\
0 & 0 & 0 & \hexstar & ... & \hexstar & \hexstar & \hexstar & \hexstar\\
0 & \hexstar & \hexstar & \hexstar & ... & \hexstar & \hexstar & \hexstar & \hexstar\\
0 & 0 & \hexstar & \hexstar & ... & \hexstar & \hexstar & \hexstar & \hexstar
\end{pmatrix}.
\]

In addition, if $\psi\in AHHol_{d}^{LF}(\mathbb{C}P^{3})\subset HHol_{d}^{LF}(\mathbb{C}P^{3})$
then by the action of an appropriate element $g\in PSp(2)$ it can
be reduced to the curve $g\psi$ with coefficient matrix of form
\[
\begin{pmatrix}\hexstar & \hexstar & \hexstar & \hexstar & ... & \hexstar & 0 & 0 & 0\\
0 & 0 & 0 & \hexstar & ... & \hexstar & \hexstar & \hexstar & \hexstar\\
0 & \hexstar & \hexstar & \hexstar & ... & \hexstar & \hexstar & 0 & 0\\
0 & 0 & \hexstar & \hexstar & ... & \hexstar & \hexstar & \hexstar & 0
\end{pmatrix}.
\]
\end{lemma}
\begin{proof}
First, the transitivity of the action of $PSp(2)$ on $\mathbb{C}P^{3}$
allows us to change the first column to the announced form. Next,
we notice that the elements of $Sp(2)$ that preserve the vector $(1,0,0,0)^{t}$
have the form
\[
\begin{pmatrix}1 & 0 & 0 & 0\\
0 & 1 & 0 & 0\\
0 & 0 & v_{1} & -\bar{v}_{2}\\
0 & 0 & v_{2} & \bar{v}_{1}
\end{pmatrix},
\]

\noindent and, since we have a matrix from $U(2)$ in the right bottom
corner, we can eliminate the last term of the second column. Other
two zeros in the second row is due to the horizontality condition
\ref{eq:2.3}. Finally, we get zeros in the right side of the coefficient
matrix due to Remark \prettyref{rem:In-the-case}.
\end{proof}

\section{\label{sec:6}Spaces of linearly full harmonic maps of degree less
than 6}

In previous sections we have shown that there is no linearly full
harmonic map of $\mathbb{R}P^{2}$ to $S^{4}$ of even degree. Also
it is clear (due to Theorem \ref{(J.-Barbosa-).}) that there are
no linearly full maps of degree less than 3. Let us describe the spaces
of such maps of degree 3 and 5. Here we have the following results.
\begin{theorem}
\label{thm:6.1}Every linearly full horizontal holomorphic curve $\psi\in AHHol_{3}^{LF}(\mathbb{C}P^{3})$
can be reduced by the action of appropriate elements $g\in PSp(2)$
and $\omega\in PSU(2)$ to the canonical form $\Psi^{3}(z)=g\psi\omega(z)=[1,-z^{3},\sqrt{3}z,\sqrt{3}z^{2}]$.
\end{theorem}
\begin{proof}
Using \prettyref{lem:5.4} and remark \prettyref{rem:In-the-case}
we get the coefficient matrix

\[
\begin{pmatrix}1 & 0 & 0 & 0\\
0 & 0 & 0 & 1\\
0 & a & 0 & 0\\
0 & 0 & -\bar{a} & 0
\end{pmatrix}.
\]

\noindent Then, acting on this curve with $\omega(z)=-z,\,g=diag\{1,1,\mathrm{exp}(-i\mathrm{arg(a)}),\mathrm{exp}(i\mathrm{arg(a)})\}$,
we get 

\[
\begin{pmatrix}1 & 0 & 0 & 0\\
0 & 0 & 0 & -1\\
0 & -a & 0 & 0\\
0 & 0 & -a & 0
\end{pmatrix},
\]

\noindent with $a\in\mathbb{R}$. Checking horizontality condition
\ref{eq:2.3} we conclude $a^{2}=3$. Thus we have two cases

\noindent 
\[
\begin{pmatrix}1 & 0 & 0 & 0\\
0 & 0 & 0 & -1\\
0 & \sqrt{3} & 0 & 0\\
0 & 0 & \sqrt{3} & 0
\end{pmatrix},\;\begin{pmatrix}1 & 0 & 0 & 0\\
0 & 0 & 0 & -1\\
0 & -\sqrt{3} & 0 & 0\\
0 & 0 & -\sqrt{3} & 0
\end{pmatrix}
\]

\noindent which can be interchanged by $g=diag\{1,1,-1,-1\}$. This
completes the proof.
\end{proof}

\begin{corollary}
\label{cor:6.2}Spaces $AHHol_{3}^{LF}(\mathbb{C}P^{3})$ and $AHarm_{3}^{\pm}(S^{4})$
are path-connected.
\end{corollary}
\begin{proof}
follows from \prettyref{thm:6.1}, \prettyref{lem:4.1} and path-connectedness
and continuity of actions of $PSp(2)$ and $PSU(2)$.
\end{proof}

Thus, the spaces of linearly full harmonic maps of $\mathbb{R}P^{2}$
to $S^{4}$ of area $6\pi$ in the induced metric have two path-connected
components $AHarm_{3}^{+}(S^{4})$ and $AHarm_{3}^{-}(S^{4})$, which
are homeomorphic through the antipodal involution of $S^{4}$. It
is worth mention that all induced metrics by these maps are canonical
since $PSp(2)$ and $PSU(2)$ act by isometries on the respective
Riemannian manifolds.
\begin{theorem}
\label{thm:6.3}Every linearly full horizontal holomorphic curve $\psi\in AHHol_{5}^{LF}(\mathbb{C}P^{3})$
by the action of appropriate elements $g\in PSp(2)$ and $\omega\in PSU(2)$
can be reduced to the canonical form $\Psi_{\eta}^{5}(z)=g\psi\omega(z)=[1+qz,-qz^{4}+z^{5},\mu z+\eta z^{2},\eta z^{3}-\mu z^{4}]$,
where $q=\frac{\eta}{\sqrt{3}}\sqrt{\frac{\eta^{2}+5}{\eta^{2}+4}},\,\mu=\frac{2}{\sqrt{3}}\sqrt{\frac{\eta^{2}+5}{\eta^{2}+4}},\,{\eta\in[0,+\infty)}$.
\end{theorem}
\begin{proof}
First, we can conclude from \prettyref{eq:3.2} that such curves have
four higher singularities. Certainly, they are two pairs of antipodal
points. We choose one of these points and rotate it to 0 by an appropriate
element of $PSU(2)$. Next, using \prettyref{lem:5.4} and Remark
\prettyref{rem:In-the-case} we obtain the coefficient matrix

\[
\begin{pmatrix}1 & q & p & 0 & 0 & 0\\
0 & 0 & 0 & \bar{p} & -\bar{q} & 1\\
0 & \mu & \eta & \nu & 0 & 0\\
0 & 0 & -\bar{\nu} & \bar{\eta} & -\bar{\mu} & 0
\end{pmatrix}.
\]

\noindent Using horizontality condition \ref{eq:2.3}, we get a
system of equations

\begin{equation}
\begin{cases}
-3p+\nu\bar{\mu}=0,\\
2q-p\bar{q}-\eta\bar{\mu}=0,\\
-3q\bar{q}+p\bar{p}+5-3\mu\bar{\mu}+\eta\bar{\eta}+\nu\bar{\nu}=0.
\end{cases}\label{eq:6.1}
\end{equation}

\noindent Using \prettyref{eq:3.1} and taking into account \ref{eq:6.1},
we get the equation of higher singularities
\begin{equation}
pz^{4}+2qz^{3}+(5-q\bar{q}-\mu\bar{\mu})z^{2}-2\bar{q}z+\frac{\mu\bar{\nu}}{3}=0.
\end{equation}

\noindent Since we have higher singularity at 0, we get $\mu\bar{\nu}=0$.
Then the first equation of \ref{thm:6.1} implies follows $p=0$.
Let us suppose that $\mu=0$, then the second equation of \ref{eq:6.1}
implies $q=0$, but this contradicts last equation of \ref{eq:6.1}.
Thus, we have $\nu=0$.

\noindent The next step is transforming non-zero coefficients into
the reals. Let $q=\tilde{q}e^{ix}$, then using multiplication by
$e^{i\nicefrac{5}{2}x}$ and actions of $\omega(z)=\mathrm{exp}(-ix)z,$\linebreak{}
 $g=diag\{e^{-\nicefrac{5}{2}x},e^{\nicefrac{5}{2}x},1,1\}$ we get
$q\in\mathbb{R}$. Further, we use the action of\linebreak{}
 $g=diag\{1,1,e^{-i\mathrm{Arg(\mu)}},e^{i\mathrm{Arg(\mu)}}\}$ and
get $\mu\in\mathbb{R}$. We conclude from the second equation of \ref{eq:6.1}
that also $\eta\in\mathbb{R}$. So, we have the coefficient matrix

\[
\begin{pmatrix}1 & q & 0 & 0 & 0 & 0\\
0 & 0 & 0 & 0 & -q & 1\\
0 & \mu & \eta & 0 & 0 & 0\\
0 & 0 & 0 & \eta & -\mu & 0
\end{pmatrix}
\]

\noindent and equations

\begin{equation}
\begin{cases}
2q-\eta\mu=0,\\
-3q^{2}+5-3\mu^{2}+\eta^{2}=0,
\end{cases}\label{eq:6.3}
\end{equation}

\noindent with real coefficients. The solutions of these equations
are 
\[
\mu=\pm\frac{2}{\sqrt{3}}\sqrt{\frac{\eta^{2}+5}{\eta^{2}+4}},\,q=\pm\frac{\eta}{\sqrt{3}}\sqrt{\frac{\eta^{2}+5}{\eta^{2}+4}},\,\eta\in\mathbb{R}.
\]
This pair of families of the solutions can be interchanged by multiplication
on $-i$ and actions of $\omega(z)=-z,\,g=diag\{i,-i,i,-i\}$. The
last step is to notice that multiplication by $-i$ and actions of
$\omega(z)=-z,\,g=diag\{i,-i,-i,i\}$ send $\eta,q,\mu$ to $-\eta,-q,\mu$,
respectively. This completes the proof.
\end{proof}

\begin{corollary}
\label{cor:6.4}The spaces $AHHol_{5}^{LF}(\mathbb{C}P^{3})$ and
$AHarm_{5}^{\pm}(S^{4})$ are path-connected.
\end{corollary}
\begin{proof}
We can see that the canonical form $\Psi_{\eta}^{5}(z)$ depends continuously
on $\eta$. Thus, there is a path connecting any canonical form $\Psi_{\eta}^{5}(z)$
with $\Psi_{0}^{5}(z)$. Then our statement follows from \prettyref{thm:6.3},
\prettyref{lem:4.1} and path-connectedness and continuity of actions
of $PSp(2)$ and $PSU(2)$.
\end{proof}

Hence the space of linearly full harmonic maps of $\mathbb{R}P^{2}$
to $S^{4}$ of area $10\pi$ in the induced metric has two path-connected
components $AHarm_{5}^{+}(S^{4})$ and $AHarm_{5}^{-}(S^{4})$, which
are homeomorphic through the antipodal involution of $S^{4}$.
\begin{remark}
\label{rem:6.5}Under the conditions of \prettyref{thm:6.3} the coefficient
matrix of the canonical form for $\eta=0$ is 
\[
\begin{pmatrix}1 & 0 & 0 & 0 & 0 & 0\\
0 & 0 & 0 & 0 & 0 & 1\\
0 & -\sqrt{\nicefrac{5}{3}} & 0 & 0 & 0 & 0\\
0 & 0 & 0 & 0 & \sqrt{\nicefrac{5}{3}} & 0
\end{pmatrix},
\]

\noindent and herewith four umbilic points with $r_{1}(p)=1$ are
glued in two umbilic points with $r_{1}(0)=2,\,r_{1}(\infty)=2$.

\noindent On the other hand if we put $\eta\rightarrow+\infty$ then
we get the following coefficient matrix
\[
\begin{pmatrix}0 & 1 & 0 & 0 & 0 & 0\\
0 & 0 & 0 & 0 & -1 & 0\\
0 & 0 & \sqrt{3} & 0 & 0 & 0\\
0 & 0 & 0 & \sqrt{3} & 0 & 0
\end{pmatrix}.
\]

\selectlanguage{english}%
\noindent This matrix gives us a map from $\mathbb{C}P^{1}=S^{2}$
to $\mathbb{C}P^{3}$ which is not defined in $z=0$, and coincides
with \foreignlanguage{american}{the canonical curve of degree 3 from
\prettyref{thm:6.1} at all other points. It means that we can define
this limit map by continuity at the point $z=0$, and this map will
be the same as the canonical curve of degree 3 with coefficient matrix
\[
\begin{pmatrix}1 & 0 & 0 & 0\\
0 & 0 & 0 & -1\\
0 & \sqrt{3} & 0 & 0\\
0 & 0 & \sqrt{3} & 0
\end{pmatrix}.
\]
It looks like we have here an annihilation of pairs of umbilic points
at $0$ and $\infty$. We look closer on this further.}

It is also worth mention here that $\eta$ parameterizes the angle
between diameters connecting pairs of antipodal singularities on $S^{2}$,
where $\eta=0$ corresponds to the angle equal to $0$, and $\eta=+\infty$
to the angle equal to $\pi$, and as we seen before this two cases
are different. We see that curves in \prettyref{thm:6.3} with different
values of $\eta$ belong to different orbits of the action of \foreignlanguage{american}{$PSp(2)$
and $PSU(2)$} since the action of $PSp(2)$
does not change higher singularities and $PSU(2)$ does not change
the angle between any pair of diameters.
\end{remark}
One can ask ``What about area? Why it discontinuously decreases?''
The answer for this question is that we have a bubbling phenomenon.
We will find the weak limit of the sequence of conformal factors of
metrics induced on $\mathbb{R}P^{2}$ to see it. In fact, the bubble
occur right at the converging point of higher singularities. For this
reason we need a more appropriate parameterization in which the pair
of higher singularities is real, symmetric with respect to zero and
converges to zero. Same calculations as in the proof of \prettyref{thm:6.3}
give us following 1-parametric family
\[
\Psi_{m}^{5}(z)=[1+pz^{2},pz^{3}+z^{5},mz+nz^{3},-nz^{2}-mz^{4}],
\]

\noindent with $n=3\sqrt{\frac{3m^{2}-5}{m^{2}+9}}$ and $p=\frac{nm}{3}=m\sqrt{\frac{3m^{2}-5}{m^{2}+9}}$,
$|m|\geq\sqrt{\frac{5}{3}}$.

Let us look at the limit when $m$ tends to $+\infty$. In this parameterization
\[
\underset{m\rightarrow+\infty}{lim}\Psi_{m}^{5}(z)=[z,z^{2},\frac{1}{\sqrt{3}},-\frac{1}{\sqrt{3}}z^{3}]
\]
is of the third degree (in the sense of Remark \prettyref{rem:6.5}). 

It is well known fact that there is just one conformal class of metrics
on $S^{2}$ and on $\mathbb{R}P^{2}$. This fact allows us to compute
a conformal factor $h(z)$ of each induced metric with respect to
the canonical metric on $S^{2}$ of area $4\pi$ ($g_{ind}=h(z)g_{can}$).
Let us apply both induced and canonical metrics to the reification
$X$ of the vector field $\sqrt{2}\frac{\partial}{\partial z}$ and
find 
\begin{equation}
h(z)=\frac{g_{ind}(X,X)}{g_{can}(X,X)}.\label{eq:6.4}
\end{equation}
The canonical metric on $S^{2}$ in terms of the isothermal coordinate
$z=x+iy$ is $g_{can}=4\frac{dx^{2}+dy^{2}}{(1+|z|^{2})^{2}}$. Hence,
$g_{can}(X,X)=\frac{4}{(1+|z|^{2})^{2}}$. We endowed $\mathbb{C}P^{3}$
with the Fubini-Study metric 
\[
g_{FS}=4\frac{|\mathbf{Z}|^{2}|d\mathbf{Z}|^{2}-(\bar{\mathbf{Z}},d\mathbf{Z})(\mathbf{Z},d\bar{\mathbf{Z}})}{|\mathbf{Z}|^{4}}=4\frac{Z_{\alpha}\bar{Z^{\alpha}}dZ_{\beta}d\bar{Z^{\beta}}-\bar{Z^{\alpha}}Z_{\beta}dZ_{\alpha}d\bar{Z^{\beta}}}{(Z_{\alpha}\bar{Z^{\alpha}})^{2}}
\]
 which has a reification $\tilde{g}_{FS}=Re\,g_{FS}$. It is easy
to check that 
\begin{equation}
\tilde{g}_{FS}(Re\,\Xi,Re\,\Xi)=\frac{1}{2}g_{FS}(\Xi,\Xi)\label{eq:6.5}
\end{equation}
for any holomorphic vector field $\Xi$. It suffices now to compute
$\tilde{g}_{FS}((\Psi_{m}^{5})_{\ast}X,(\Psi_{m}^{5})_{\ast}X)=g_{ind;m}(X,X)$
for finding the denominator of \prettyref{eq:6.4}. Indeed, twistor
fibration is a Riemannian submersion and $\Psi_{m}^{5}$ is orthogonal
to the fibers as it mentioned before. We have 
\[
\tilde{g}_{FS}((\Psi_{m}^{5})_{\ast}X,(\Psi_{m}^{5})_{\ast}X)=4\frac{|\Psi_{m}^{5}|^{2}|\frac{\partial}{\partial z}\Psi_{m}^{5}|^{2}-|(\overline{\Psi_{m}^{5}},\frac{\partial}{\partial z}\Psi_{m}^{5})|^{2}}{|\Psi_{m}^{5}|^{4}}=
\]
\[
=4\frac{\partial}{\partial\bar{z}}\frac{\partial}{\partial z}\ln|\Psi_{m}^{5}|^{2}=\triangle\ln|\Psi_{m}^{5}|^{2}
\]
 according to \prettyref{eq:6.5}. We do further calculations in the
polar coordinates $z=r(\cos\varphi+i\sin\varphi)$, $\triangle=\frac{1}{r}\frac{\partial}{\partial r}(r\frac{\partial}{\partial r})+\frac{1}{r^{2}}\frac{\partial^{2}}{\partial\varphi^{2}}$.
Let $h_{m}(z)$ be the conformal factor for the metric induced by
$\Psi_{m}^{5}$.
\begin{proposition}
\label{Prop. 3}Conformal factors $h_{m}(z)$ converge as $m\rightarrow+\infty$
to the $4\pi\delta(x,y)+3$ in the weak sense ($\delta(x,y)$ is Dirac
delta). 
\end{proposition}
\begin{proof}
First, we find out that $|\Psi_{m}^{5}|^{2}=(1+r^{2})^{3}w_{m}(r,\varphi)$
with 
\[
w_{m}(r,\varphi)=(2pr^{2}\cos2\varphi+r^{4}+(m^{2}-3)r^{2}+1).
\]
 Then, clearly, 
\[
\triangle\ln|\Psi_{m}^{5}|^{2}=\triangle\ln w_{m}(r,\varphi)+3\triangle\ln(1+r^{2})=\triangle\ln w_{m}(r,\varphi)+\frac{12}{(1+r^{2})^{2}}
\]
 and 
\[
h_{m}(z)=\frac{(1+r^{2})^{2}}{4}\triangle\ln w_{m}(r,\varphi)+3.
\]
 Thus, we have to show that $\tilde{h}_{m}(r,\varphi)=\frac{(1+r^{2})^{2}}{4}\triangle\ln w_{m}(r,\varphi)$
converges weakly to $4\pi\delta(x,y)$. Since $\intop_{0}^{2\pi}d\varphi\intop_{0}^{1}h_{m}(z)\frac{4rdr}{(1+r^{2})^{2}}=10\pi$
is an area of $\mathbb{R}P^{2}$ it sufficient to show that\linebreak{}
 
\[
\underset{m\rightarrow+\infty}{\lim}\intop_{0}^{2\pi}d\varphi\intop_{\frac{1}{m}}^{1}\tilde{h}_{m}(r,\varphi)\frac{4rdr}{(1+r^{2})^{2}}=\underset{m\rightarrow+\infty}{\lim}\intop_{0}^{2\pi}d\varphi\intop_{\frac{1}{m}}^{1}r\triangle\ln w_{m}(r,\varphi)dr=0.
\]
 First, we note that $\intop_{0}^{2\pi}\frac{\partial^{2}}{\partial\varphi^{2}}\ln w_{m}(r,\varphi)d\varphi=\frac{\partial}{\partial\varphi}\ln w_{m}(r,\varphi)|_{0}^{2\pi}=0$.
Next, 
\[
\tilde{w}_{m}(r,\varphi)=r\frac{\partial}{\partial r}\ln w_{m}(r,\varphi)=2+2\frac{r^{4}-1}{w_{m}(r,\varphi)},
\]
\[
\tilde{w}_{m}(0,\varphi)=0,
\]

\[
\tilde{w}_{m}(1,\varphi)-\tilde{w}_{m}(\frac{1}{m},\varphi)=\frac{1}{p}(1-\frac{1}{m^{2}})\frac{1}{\cos2\varphi+\frac{1}{2pm^{2}}+\frac{m^{2}-2}{2p}}<\frac{2}{m^{2}-2p-2},
\]
 and thus 
\[
\intop_{0}^{2\pi}d\varphi\intop_{\frac{1}{m}}^{1}r\frac{1}{r}\frac{\partial}{\partial r}(r\frac{\partial}{\partial r})\ln w_{m}(r,\varphi)dr=\intop_{0}^{2\pi}d\varphi\left(r\frac{\partial}{\partial r}\ln w_{m}(r,\varphi)\right)|_{\frac{1}{m}}^{1}=
\]
\[
=\intop_{0}^{2\pi}d\varphi\left(\tilde{w}_{m}(1,\varphi)-\tilde{w}_{m}(\frac{1}{m},\varphi)\right)<\intop_{0}^{2\pi}\frac{2}{m^{2}-2p-2}d\varphi\xrightarrow[m\rightarrow+\infty]{}0.
\]
\end{proof}

This means that we have here the bubbling phenomenon, i.e. the ``bubble''
of area $4\pi$ is growing at the point $0$ whilst the metric on
$\mathbb{R}P^{2}$ converges to canonical metric of area $6\pi$.
The word ``bubble'' means $S^{2}$ with canonical metric. And actually,
we observe that the area is preserved. We see here that ``bubble''
appear right at the point $0$ where we defined the limit map by continuity\foreignlanguage{american}{
in }\prettyref{rem:6.5}\foreignlanguage{american}{. In fact, pairs
of higher singularities can annihilate in any point (for example,
we can just send $0$ to other point by the action of $PSU(2)$).
This means that there are other copies of $AHHol_{3}^{LF}(\mathbb{C}P^{3})$
in the closure of $AHHol_{5}^{LF}(\mathbb{C}P^{3})$ in the space
$P((\mathbb{C}[z]_{5})^{4})$ and they differ by the point where ``bubble''
appears, i.e. regular parts of conformal factors of induced metrics
are the same but singular parts are different.}

\selectlanguage{american}%
Let us now look at the closure of $AHHol_{d}^{LF}(\mathbb{C}P^{3})$
in the space $P((\mathbb{C}[z]_{d})^{4})$. We call the relative complement
of $AHHol_{d}^{LF}(\mathbb{C}P^{3})$ in the closure of $AHHol_{d}^{LF}(\mathbb{C}P^{3})$
in $P((\mathbb{C}[z]_{d})^{4})$ the \textit{boundary} of $AHHol_{d}^{LF}(\mathbb{C}P^{3})$
in the space $P((\mathbb{C}[z]_{d})^{4})$. We have the following
lemma for the boundary of $AHHol_{d}^{LF}(\mathbb{C}P^{3})$.
\begin{lemma}
A point in the boundary of $AHHol_{d}^{LF}(\mathbb{C}P^{3})$ could
be identified with a linearly full invariant horizontal holomorphic
curve of degree less than d.
\end{lemma}
\begin{proof}
Recall that $AHHol_{d}^{LF}(\mathbb{C}P^{3})$ is the subspace in
the space $P((\mathbb{C}[z]_{d})^{4})$ defined by the number of conditions.
Let us figure out which of these conditions are held by the closure
points of $AHHol_{d}^{LF}(\mathbb{C}P^{3})$. Horizontality and invariance
conditions are held since they are defined by algebraic equalities.
It is also easy to show that if a limiting curve (a point of the closure)
are not linearly full then all the coefficients of such curve must
vanish. Hence, the closure points represent linearly full maps. We
see that there are two conditions which could not be preserved by
closure points, - the condition of degree, and the condition of coprime
polynomials. For polynomials that are not coprime we can define the
required curve by dividing all these polynomials by greatest common
divisor, which is the same as to define the limit map by continuity
at the points where greatest common divisor vanishes.
\end{proof}

\selectlanguage{english}%
\begin{remark}
It is now fair to assume that ``bubbles'' occur at the zeros of
greatest common divisor, i.e. we have here the same phenomenon as
in Proposition \ref{Prop. 3}.
\end{remark}
Now we need to define the space \foreignlanguage{american}{$AHHol_{\leq5}^{LF}(\mathbb{C}P^{3})$
of linearly full invariant horizontal holomorphic curves of degree
less or equal to 5. We define it as the closure of $AHHol_{5}^{LF}(\mathbb{C}P^{3})$
in the space $P((\mathbb{C}[z]_{5})^{4})$.}
\selectlanguage{american}%
\begin{proposition}
\label{prop.4}The boundary of $AHHol_{5}^{LF}(\mathbb{C}P^{3})$
in the space $P((\mathbb{C}[z]_{5})^{4})$ is the orbit of the point
with $\eta=+\infty$ in Remark \prettyref{rem:6.5} under the action
of $PSp(2)$ and $PSU(2)$.
\end{proposition}
\selectlanguage{english}%
\begin{proof}
Let us denote by $L\subset AHHol_{5}^{LF}(\mathbb{C}P^{3})$ the subspace
of canonical forms from \prettyref{thm:6.3}, i.e. the subspace consisting
of points $\Psi_{\eta}^{5}\in AHHol_{5}^{LF}(\mathbb{C}P^{3})$, with
$\eta\in[0,+\infty)$. Then the only limit point of $L$ in $P((\mathbb{C}[z]_{5})^{4})$
which does not lie in $L$ is that \foreignlanguage{american}{with
$\eta=+\infty$ in Remark \prettyref{rem:6.5}.} The groups \foreignlanguage{american}{$PSp(2)$
and $PSU(2)$ act on the larger space $P((\mathbb{C}[z]_{5})^{4})\supset AHHol_{5}^{LF}(\mathbb{C}P^{3})$
and this action is continuous, in the sense that the action map
\[
\mathscr{A}:\:PSp(2)\times P((\mathbb{C}[z]_{5})^{4})\times PSU(2)\rightarrow P((\mathbb{C}[z]_{5})^{4}),\;\mathscr{A}:\:(g,\psi,\omega)\mapsto g\psi\omega
\]
is continuous. Let us look at a boundary point $\psi$. It is the
limit of the sequence $\{\psi_{n}\}\subset AHHol_{5}^{LF}(\mathbb{C}P^{3})$.
By the \prettyref{thm:6.3} we have sequences $\{(g_{n},\omega_{n})\}\subset PSp(2)\times PSU(2)$
and $\{\Psi_{\eta_{n}}^{5}\}\subset L\subset AHHol_{5}^{LF}(\mathbb{C}P^{3})$,
such that $\psi_{n}=g_{n}\Psi_{\eta_{n}}^{5}\omega_{n}$. Since $PSp(2)\times P((\mathbb{C}[z]_{5})^{4})\times PSU(2)$
is compact and metrizable, we have a subsequence $\{(g_{k},\Psi_{\eta_{k}}^{5},\omega_{k})\}$,
which converges to a point $(g,\Psi_{\infty}^{5},\omega)$, with $\Psi_{\infty}^{5}=\underset{k\rightarrow+\infty}{\lim}\Psi_{\eta_{k}}^{5}$.
Note that the subsequence $\{\psi_{k}=g_{k}\Psi_{\eta_{k}}^{5}\omega_{k}\}$
converges to the point $\psi$. Since $\mathscr{A}$ is continuous,
we have $\mathscr{A}(g,\Psi_{\infty}^{5},\omega)=g\psi\omega$. Clearly,
the point $\Psi_{\infty}^{5}\in P((\mathbb{C}[z]_{5})^{4})$ does
not lie in $L$, but it is a limit point of $L$. This completes the
proof of the inclusion $boundary\subset orbit$. The proof of the
inverse inclusion is trivial.}
\end{proof}

\selectlanguage{american}%
\begin{theorem}
\label{Thm6.8}The space $AHHol_{\leq5}^{LF}(\mathbb{C}P^{3})$ is
path-connected.
\end{theorem}
\begin{proof}
follows immediately from the definition of $AHHol_{\leq5}^{LF}(\mathbb{C}P^{3})$,
Corollary \prettyref{cor:6.4}, Proposition \ref{prop.4}, and path-connectedness
of the closure of $L\subset P((\mathbb{C}[z]_{5})^{4})$.
\end{proof}

\section{\label{sec:7}Linearly full harmonic maps of the two-dimensional
sphere to the four-dimensional sphere of degree less than 6.}

In this section we summarize known results for linearly full harmonic
maps $S^{2}\rightarrow S^{4}$ of degree less than 6. In this case
we have groups $PSp(2,\mathbb{C})$ and $PSL(2,\mathbb{C})$ acting
on the spaces $HHol_{d}^{LF}(\mathbb{C}P^{3})$ by postcomposition
and precomposition respectively. Note that these groups are path-connected
\citep{key-10}. For the beginning we provide a useful lemma from
\citep{key-6}.
\begin{lemma}
\textup{\citep{key-6}}\label{lem:7.1} The group $PSp(2,\mathbb{C})$
acts transitively on $\mathbb{C}P^{3}$ and the stabilizer of the
point $[1,0,0,0]^{t}$ consists of the matrices of the form

\[
\begin{pmatrix}\xi & \eta & \xi(\gamma\lambda-\alpha\mu) & \xi(\delta\lambda-\beta\mu)\\
0 & \zeta & 0 & 0\\
0 & \lambda & \alpha & \beta\\
0 & \mu & \gamma & \delta
\end{pmatrix},\quad\alpha\delta-\beta\gamma=\xi\zeta=1.
\]
\end{lemma}
It is clear that there are no harmonic maps of degree less than 3.
So we start with degree 3. Due to \prettyref{eq:3.2} we deduce that
such maps can not have higher singularities. In this simple case we
have following result.
\begin{theorem}
\label{thm:7.2}Every linearly full horizontal holomorphic curve $\psi\in HHol_{3}^{LF}(\mathbb{C}P^{3})$
can be reduced by the action of appropriate elements $g\in PSp(2,\mathbb{C})$
and $\omega\in PSL(2,\mathbb{C})$ to the canonical form $\tilde{\Psi}^{3}(z)=g\psi\omega(z)=[1,-z^{3},\sqrt{3}z,\sqrt{3}z^{2}]$.
\end{theorem}
\begin{proof}
Using \prettyref{lem:5.4} we get coefficient matrix of form

\begin{equation}
\begin{pmatrix}1 & a & b & c\\
0 & 0 & 0 & d\\
0 & r & p & q\\
0 & 0 & s & t
\end{pmatrix}.
\end{equation}

\noindent Next, we use \prettyref{lem:7.1} several times. After the
each time, we rename all non-zero coefficients as they were before.
So, our sequence of transformation matrices from \prettyref{lem:7.1}
is following

\begin{equation}
\begin{array}{c}
\begin{pmatrix}1 & 0 & \nicefrac{t}{d} & 0\\
0 & 1 & 0 & 0\\
0 & 0 & d^{-1} & 0\\
0 & -t & 0 & d
\end{pmatrix},\,\begin{pmatrix}1 & 0 & 0 & -\nicefrac{q}{d}\\
0 & 1 & 0 & 0\\
0 & -q & d & 0\\
0 & 0 & 0 & d^{-1}
\end{pmatrix},\,\begin{pmatrix}1 & 0 & 0 & 0\\
0 & 1 & 0 & 0\\
0 & 0 & s & -p\\
0 & 0 & 0 & s^{-1}
\end{pmatrix},\\
\begin{pmatrix}d & -c & 0 & 0\\
0 & d^{-1} & 0 & 0\\
0 & 0 & 1 & 0\\
0 & 0 & 0 & 1
\end{pmatrix},\,\begin{pmatrix}\sqrt{-d} & 0 & 0 & 0\\
0 & \sqrt{-d}^{-1} & 0 & 0\\
0 & 0 & \sqrt{\frac{s}{r}} & 0\\
0 & 0 & 0 & \sqrt{\frac{r}{s}}
\end{pmatrix}.
\end{array}\label{eq:7.2}
\end{equation}

\noindent This sequence of transformations gives us following sequence
of coefficient matrices

\begin{equation}
\begin{array}{c}
\begin{pmatrix}1 & a & b & c\\
0 & 0 & 0 & d\\
0 & r & p & 0\\
0 & 0 & s & t
\end{pmatrix}\rightarrow\begin{pmatrix}1 & a & b & c\\
0 & 0 & 0 & d\\
0 & r & p & 0\\
0 & 0 & s & 0
\end{pmatrix}\rightarrow\begin{pmatrix}1 & a & b & c\\
0 & 0 & 0 & d\\
0 & r & 0 & 0\\
0 & 0 & s & 0
\end{pmatrix}\rightarrow\\
\rightarrow\begin{pmatrix}1 & a & b & 0\\
0 & 0 & 0 & d\\
0 & r & 0 & 0\\
0 & 0 & s & 0
\end{pmatrix}\rightarrow\begin{pmatrix}1 & a & b & 0\\
0 & 0 & 0 & -1\\
0 & s & 0 & 0\\
0 & 0 & s & 0
\end{pmatrix}.
\end{array}
\end{equation}

\noindent Now we use the horizontality condition \prettyref{eq:2.3}
and conclude that $a=b=0,\,\frac{s^{2}}{3}=1$. We have $s=\pm\sqrt{3}$
and these $\pm$ cases can be interchanged by $g=diag\{1,1,-1,-1\}$.
This completes the proof.
\end{proof}

We have the following.
\begin{corollary}
\label{cor7.3}Spaces $HHol_{3}^{LF}(\mathbb{C}P^{3})$ and $Harm_{3}^{\pm}(S^{4})$
are path-connected.
\end{corollary}
\begin{proof}
follows from \prettyref{thm:7.2}, \prettyref{lem:4.1} and path-connectedness
and continuity of actions of $PSp(2,\mathbb{C})$ and $PSL(2,\mathbb{C})$.
\end{proof}

Here we see that the canonical form $\tilde{\Psi}^{3}(z)$ for $S^{2}$
is equal to the canonical form $\Psi^{3}(z)$ for $\mathbb{R}P^{2}$.

We continue with the case of degree 4. We conclude from \prettyref{eq:3.2}
that such maps can have either 1 or 2 higher singularities. We will
see later that case of one higher singularity is impossible. In fact
we have following.
\begin{theorem}
\label{thm:7.4}Every linearly full horizontal holomorphic curve $\psi\in HHol_{4}^{LF}(\mathbb{C}P^{3})$
can be reduced by the action of appropriate elements $g\in PSp(2,\mathbb{C})$
and $\omega\in PSL(2,\mathbb{C})$ to the canonical form 
\[
\tilde{\Psi}_{a}^{4}(z)=g\psi\omega(z)=[1+az,-z^{4},\sqrt{2}z(1+\nicefrac{3}{2}az),\sqrt{2}z^{3}],\,a\in\mathbb{R}.
\]
\end{theorem}
\begin{proof}
We move one of the higher singularities to zero using an appropriate
element $\omega\in PSU(2)$. Next we use \prettyref{lem:5.4} and
get the coefficient matrix of form

\begin{equation}
\begin{pmatrix}1 & a & b & e & c\\
0 & 0 & 0 & f & d\\
0 & r & m & p & q\\
0 & 0 & l & s & t
\end{pmatrix}.
\end{equation}

\noindent We conclude that $rl=0$ because of \prettyref{eq:3.1}
and higher singularity at zero. It is clear that without loss of generality
we may assume that $l=0$. Horizontality condition \prettyref{eq:2.3}
gives us $f=-\frac{rl}{3}=0$. So, in fact, our coefficient matrix
is of the form

\begin{equation}
\begin{pmatrix}1 & a & b & e & c\\
0 & 0 & 0 & 0 & d\\
0 & r & m & p & q\\
0 & 0 & 0 & s & t
\end{pmatrix}.
\end{equation}

\noindent Next we use sequence of transformations \prettyref{eq:7.2}
and get the coefficient matrix

\begin{equation}
\begin{pmatrix}1 & a & b & e & 0\\
0 & 0 & 0 & 0 & -1\\
0 & s & m & 0 & 0\\
0 & 0 & 0 & s & 0
\end{pmatrix}.
\end{equation}

\noindent For this matrix horizontality condition \prettyref{eq:2.3}
implies $b=e=0,\,3a=ms,\,s^{2}=2.$ So $s=\pm\sqrt{2}$ and these
two cases can be interchanged by $g=diag\{1,1,-1,-1\}$. We choose
$s=\sqrt{2}$, then $m=\nicefrac{3}{\sqrt{2}}a.$ Now let $x=\mathrm{arg(\mathit{a})}$,
then using multiplication by $e^{i2x}$ and actions of $g=\{e^{-i2x},e^{i2x},e^{-ix},e^{ix}\},\,\omega(z)=e^{-ix}z$
we get $a\in\mathbb{R}$. This completes the proof.
\end{proof}

\begin{remark}
\label{rem:7.5}In this theorem we see that the canonical form $\tilde{\Psi}_{a}^{4}(z)$
is continuously dependent on $a$, therefore we can deform each one
to the form $\tilde{\Psi}_{0}^{4}(z)=[1,-z^{4},\sqrt{2}z,\sqrt{2}z^{3}]$.
Equation \prettyref{eq:3.1} describing the higher singularities for
these curves is $z(az+2)=0$, therefore second singularity is at the
point $z_{0}=-\frac{2}{a}$. Note that second singularity for $\tilde{\Psi}_{0}^{4}(z)$
is at infinity. On the other hand, if we want to move second singularity
to zero, we must move $a$ to infinity, but then degree of the curve
is decreasing to $0$. So we conclude that such curves can not have
one higher singularity. Therefore, in fact, we can deform every curve
of degree 4 to $\tilde{\Psi}_{0}^{4}(z)$ by some element $g\in PSp(2,\mathbb{C})$
and Moebius transformation $\omega\in PSL(2,\mathbb{C})$ which moves
higher singularities to the points $0$ and $\infty$. In fact, we
have here a bubbling phenomenon too. But an $a$-parametric family
of curves degenerates to the point where the bubble grows as ``$a$''
tends to $\infty$. In other words, the whole area of $a$-parametric
family flows off to the area of the bubble.
\end{remark}
This remark immediately leads us to the following.
\begin{corollary}
\label{Cor7.4}Spaces $HHol_{4}^{LF}(\mathbb{C}P^{3})$ and $Harm_{4}^{\pm}(S^{4})$
are path-connected.
\end{corollary}
\begin{proof}
follows from \prettyref{thm:7.4}, \prettyref{lem:4.1}, remark \prettyref{rem:7.5}
and path-connectedness and continuity of actions of $PSp(2,\mathbb{C})$
and $PSL(2,\mathbb{C})$.
\end{proof}

Note that canonical forms for cases of degrees 3 and 4 coincide with
the following theorem from \citep{key-6}.
\begin{theorem}
\textup{\citep{key-6}}\label{thm:7.7} Let $\psi:\,S^{2}\rightarrow\mathbb{C}P^{3}$
be a linearly full horizontal holomorphic curve of degree $d$ which
has at most two higher singularities. Then there exist elements $g\in PSp(2,\mathbb{C})$
and $\omega\in PSL(2,\mathbb{C})$ such that $\Psi(z)=g\psi\mu(z)=[1,k_{2}z^{2k_{1}+k_{2}},-(2k_{1}+k_{2})z^{k_{1}},z^{k_{1}+k_{2}}]$,
for some positive integers $k_{1},k_{2}$ with $2k_{1}+k_{2}=d$.
The higher singularities of $\Psi(z)$, if any, occur at points $0$
and $\infty$.
\end{theorem}
We start description of maps of degree 5 with the simple case when
$\psi$ has at most two higher singularities.
\begin{theorem}
\label{thm:7.8}Let $\psi:\,S^{2}\rightarrow\mathbb{C}P^{3}$ be a
linearly full horizontal holomorphic curve of degree $5$ which has
at most two higher singularities. Then there exist elements $g\in PSp(2,\mathbb{C})$
and $\omega\in PSL(2,\mathbb{C})$ reducing $\psi$ to the one of
the following form

\begin{equation}
\tilde{\Psi}_{1}^{5}(z)=[1,z^{5},2z^{2},-\frac{5}{2}z^{3}],\label{eq:7.7}
\end{equation}

\begin{equation}
\tilde{\Psi}_{2}^{5}(z)=[1,z^{5},z,-\frac{5}{3}z^{4}].\label{eq:7.8}
\end{equation}
\end{theorem}
\begin{proof}
From \prettyref{thm:7.7} we get two cases $\Psi(z)=[1,z^{5},-5z^{2},z^{3}]$
and $\Psi(z)=[1,3z^{5},-5z,z^{4}]$ which can be reduced to \ref{eq:7.7}
and \ref{eq:7.8} by elements\linebreak{}
 $g_{1}=diag\{1,1,-\frac{2}{5},-\frac{5}{2}\}$ and $g_{2}=diag\{\sqrt{3},\frac{1}{\sqrt{3}},-\frac{\sqrt{3}}{5},-\frac{5}{\sqrt{3}}\}$
respectively.
\end{proof}

It is worth mentioning that, in fact, both of forms from this theorem
have two higher singularities.

So, we can consider the case of at least two higher singularities.
In this most general case we have following result.
\begin{theorem}
\textup{\citep{key-6}}\label{thm:7.9} Let $\psi:\,S^{2}\rightarrow\mathbb{C}P^{3}$
be a linearly full horizontal holomorphic curve of degree $5$. Then
there exist elements $g\in PSp(2,\mathbb{C})$ and $\omega\in PSL(2,\mathbb{C})$
such that 
\begin{equation}
\tilde{\Psi}^{5}(z)=g\psi\mu(z)=[1+az,(h+z)z^{4},(r+lz)z,(m+sz)z^{3}]\label{eq:7.9}
\end{equation}
 canonical linearly full horizontal holomorphic curve, where $a,h,r,l,m,s$
satisfy 
\begin{equation}
2a+ls=0,\,5+3ah+lm+3rs=0,\,2h+rm=0.\label{eq:7.10}
\end{equation}
\end{theorem}
\begin{remark}
\label{rem:7.10}Degree of curve from parametric family \{\ref{eq:7.9},
\ref{eq:7.10}\} decreases iff $ah=1,\,rs=lm=-2,\,hl=r$ or, equivalently,
$ah=1,\,rs=lm=-2,\,ar=l$. In this case degree of the curve is 4 and
we have a bubbling phenomenon here. Two of higher singularities are
at points $0$ and $\infty$, and another two are meeting at point
$-1$ and are annihilating with the birth of bubble of area $4\pi$.
\end{remark}
\begin{lemma}
\label{lem:7.11}Every canonical form $\tilde{\Psi}^{5}(z)$ from
\prettyref{thm:7.9} can be deformed to canonical form $\tilde{\Psi}_{0}^{5}(z)=[1+z,(2+z)z^{4},(1+2z)z,-(4+z)z^{3}]$.
\end{lemma}
\begin{proof}
The first step is in showing that there exists deformation of $\tilde{\Psi}_{1}^{5}(z)$
to $\tilde{\Psi}_{2}^{5}(z)$ from \prettyref{thm:7.8}. We will do
it by deformation both of $\tilde{\Psi}_{1}^{5}(z)$ and $\tilde{\Psi}_{2}^{5}(z)$
to $\tilde{\Psi}_{0}^{5}(z)$. These deformations are given by $\tilde{\Psi}_{1}^{5}(z,t)=[1+tz,(2t+z)z^{4},(t+2z)z,(\frac{-5-3t^{2}}{2}+tz)z^{3}],\,t\in[0,1]$
and $\tilde{\Psi}_{2}^{5}(z)=[1+tz,(2t+z)z^{4},(1+2tz)z,(-4t+\frac{2t^{2}-5}{3}z)z^{3}],\,t\in[0,1]$
respectively.

\noindent The second step is in deforming every curve $\tilde{\Psi}^{5}(z)$
to one of the two curves $\tilde{\Psi}_{1}^{5}(z)$ and $\tilde{\Psi}_{2}^{5}(z)$.
Here we have three cases

1) $rs=0,\,lm=-5$ and deformation to $\tilde{\Psi}_{1}^{5}(z)$ is
given by $\tilde{\Psi}_{1}^{5}(z,t)=[1+a(1-t)z,(h(1-t)+z)z^{4},(r(1-t)+lz)z,(m+s(1-t)z)z^{3}],\,t\in[0,1]$
and then applying $g_{1}=diag\{1,1,\frac{2}{l},-\frac{5}{2m}\}$.

2) $lm=0,\,rs=-\frac{5}{3}$ and deformation to $\tilde{\Psi}_{2}^{5}(z)$
is given by $\tilde{\Psi}_{2}^{5}(z,t)=[1+a(1-t)z,(h(1-t)+z)z^{4},(r+l(1-t)z)z,(m(1-t)+sz)z^{3}],\,t\in[0,1]$
and then applying $g_{2}=diag\{1,1,\frac{1}{r},-\frac{5}{3s}\}$.

3) $r,l,m,s\neq0$ and deformation to $\tilde{\Psi}_{2}^{5}(z)$ is
given by $\tilde{\Psi}_{2}^{5}(z,t)=[1+a(1-t)z,(h(1-t)+z)z^{4},(r+l(1-t)z)z,(m(1-t)-\frac{5+(3ah+lm)(1-t)^{2}}{3r}z)z^{3}],\,t\in[0,1]$
and then applying $g_{2}=diag\{1,1,\frac{1}{r},r\}$.

\noindent It is clear from Remark \prettyref{rem:7.10} that the degree
is preserved throughout all of these deformations. This completes
the proof.
\end{proof}

\begin{corollary}
\label{Cor7.5}Spaces $HHol_{5}^{LF}(\mathbb{C}P^{3})$ and $Harm_{5}^{\pm}(S^{4})$
are path-connected.
\end{corollary}
\begin{proof}
follows from \prettyref{thm:7.9}, \prettyref{lem:4.1}, \prettyref{lem:7.11}
and path-connectedness and continuity of actions of $PSp(2,\mathbb{C})$
and $PSL(2,\mathbb{C})$.
\end{proof}

For the closures of spaces $HHol_{d}^{LF}(\mathbb{C}P^{3})$ in $P((\mathbb{C}[z]_{d})^{4})$
we could see that limit points could represent not linearly full maps.
As an example, one could just take any $\psi=[f_{1},f_{2},f_{3},f_{4}]\in HHol_{d}^{LF}(\mathbb{C}P^{3})$
and construct 1-parametric family $\psi_{t}=[f_{1},tf_{2},f_{3},tf_{4}]\in HHol_{d}^{LF}(\mathbb{C}P^{3})$
for all positive $t$, which obviously converges to not linearly full
map as $t\rightarrow0$. 

Our explicit calculations of this section agree with the works of
B. Loo \citep{key-12} and M. Kotani \citep{key-13}. Here we end
our summary of known results for harmonic maps of $S^{2}$ to $S^{4}$.
For more details one can refer to \citep{key-6}.

\textbf{Acknowledgements:} Author is thankful to Alexei V. Penskoi for attaching the author's attention to this problem, helpful remarks and useful conversations.

\noindent \textit{\medskip{}
}

\end{document}